\documentclass[nothms,english,utf8,12pt,upmath]{aclart}
\DeclareUnicodeCharacter{00A0}{}
\usepackage{macros}

\def\index#1{}\def\glossary#1{}
\begin{document}

\parindent=1em

\title [Heights, equidistribution, and the Bogomolov conjecture]{Arakelov geometry, heights, equidistribution, and the Bogomolov conjecture}
\hypersetup{pdftitle={Arakelov geometry, heights, equidistribution,
and the Bogomolov conjecture}}
 \author{Antoine Chambert-Loir}
\address{Université Paris Diderot, Sorbonne Université, CNRS, Institut de Mathématiques de Jussieu-Paris Rive Gauche, IMJ-PRG, F-75013, Paris, France}
 
 \email{Antoine.Chambert-Loir@math.univ-paris-diderot.fr}


%


\thanks{Preliminary version (29 June 2018; revised March 2019; re-revised
June 2020)
of lectures delivered at Institut Fourier, Grenoble (France)
during the summer school ``Géométrie d'Arakelov et applications diophantiennes'' (12--30 June 2017).
An updated version of this text may  be  accessible online at address
\url{https://webusers.imj-prg.fr/~antoine.chambert-loir/publications/papers/grenoble.pdf}.
}
\maketitle

\begin{flushright}
\itshape
To the memory of Lucien Szpiro
\end{flushright}

\tableofcontents

\section{Introduction}

Fermat's method of infinite descent studies the solutions to diophantine equations by constructing, from a given solution of a diophantine equation, a smaller solution, and ultimately deriving a contradiction.
In order to formalize the intuitive notion of ``size'' 
of an algebraic solution of a diophantine equation,
\cite{northcott1950} and \cite{weil1951} have introduced 
the notion of \emph{height} of an algebraic point of an algebraic
variety defined over a number field
and established their basic functorial properties,
using the decomposition theorem of \cite{weil1929}.
The \emph{height machine} is now an important tool 
in modern diophantine geometry.

The advent of arithmetic intersection theory 
with~\cite{arakelov1974} and, above all,
its extension in any dimension by~\cite{gillet-s90} (``Arakelov geometry'')
has led \cite{faltings1991} to extend the concept further by introducing
the height of a subvariety, defined in pure analogy with its degree,
replacing classical intersection theory with arithmetic intersection theory.
This point of view has been developed in great depth by~\cite{bost-g-s94}
and~\cite{zhang95}. 

Although I shall not use it in these notes,
I also mention the alternative viewpoint of~\cite{philippon91}
who defines the height of a subvariety as the height 
of the coefficients-vector of its ``Chow form''.

The viewpoint of adelic metrics introduced in~\cite{zhang95b} is strengthened by the introduction of Berkovich spaces in this context,
based on~\cite{gubler1998}, and leading to the definition 
by~\cite{chambert-loir2006} of measures at all places analogous
to product of Chern forms at the archimedean place. 

We describe these notions in the first sections of the text.
We then present the equidistribution theorem of~\cite{szpiro-u-z97}
and its extension by~\cite{yuan2008}.
Finally, we use these ideas to explain the proof of Bogomolov's conjecture,
following~\cite{ullmo98} and~\cite{zhang98}.

\medskip

I thank the editors of this volume for their invitation
to lecture on this topic in Grenoble,
and  the referees for their patient comments that 
helped me improve this survey.

\medskip

I started my graduate education under the guidance of Lucien Szpiro.
He introduced me to heights, abelian varieties and Arakelov geometry
almost 30 years ago.
As a sad coincidence, he passed away a few weeks before I was to make
the final modifications to this text, at a time when
a viral pandemy locked down half of humanity.
I dedicate this survey to his memory.

\def\ad{{\mathrm{ad}}}
\def\adm{{\mathrm{adm}}}
\def\an{{\mathrm{an}}}
\def\cpct{{\mathrm{c}}}
\def\resp{\emph{resp.}\xspace}
\section{Arithmetic intersection numbers}

\subsection{}
Let $\mathscr X$ be a proper flat scheme over~$\Z$. 
For every integer~$d\geq0$, let $\mathrm Z_d(\mathscr X)$
be the group of $d$-cycles on~$\mathscr X$: it is the free abelian group
generated by integral closed subschemes of dimension~$d$.
\begin{rema}
Let $f\colon\mathscr X\to\Spec(\Z)$ be the structural morphism.
By assumption, $f$ is proper so that the image of an integral closed subscheme~$Z$ of~$\mathscr X$ is again an integral closed subscheme of~$\Spec(\Z)$.
There are thus two cases: \begin{enumerate}
\item Either $ f(Z)=\Spec(\Z)$, in which case we say that $Z$ is horizontal;\index{Horizontal subscheme}
\item Or $f(Z)=\{ (p)\}$ for some prime number~$p$, in which case we say that $Z$ is vertical.\index{Vertical subscheme}
\end{enumerate}
\end{rema}

\subsection{}
The set $\mathscr X(\C)$ of complex points of~$\mathscr X$
has a natural structure of a complex analytic space,
smooth if and only if $\mathscr X_\Q$ is regular.
This gives rise to the notions of continuous, resp. smooth,
resp. holomorphic function on~$\mathscr X(\C)$: by definition, this is
a function which,
for every local embedding of an open subset~$U$
of~$\mathscr X(\C)$ into~$\C^n$, extends to a continuous,
resp. smooth, resp. holomorphic function around the image of~$U$.

Let $\mathscr L$ be a line bundle on~$\mathscr X$.
A \emph{hermitian metric} on~$\mathscr L$ is the datum,
for every open subset~$U$ of~$\mathscr X(\C)$ and every section
$s\in\Gamma(U,\mathscr L)$ of a continuous function $\norm s\colon U\to\R_+$,
subject to the following conditions:
\begin{enumerate}
\item For every open subset~$V$ of~$U$, one has $\norm {s|_V}=\norm{s}|_V$;
\item For every holomorphic function $f\in\mathscr O_{\mathscr X}(U)$,
one has $\norm{fs}=\abs f \norm s$;
\item If $s$ does not vanish on~$U$, then the function~$\norm s$ is
strictly positive and smooth.
\end{enumerate}
A \emph{hermitian line bundle}~$\overline{\mathscr L}$
on~$\mathscr X$ is a line bundle~$\mathscr L$
endowed with a hermitian metric.

With respect to the tensor product of underlying line bundles
and the tensor product of hermitian metrics, the set
of isomorphism classes
of hermitian line bundles on~$\mathscr X$ is an abelian group,
denoted by $\hPic(\mathscr X)$.
This group fits within an exact sequence of abelian groups:
\begin{equation}
 \Gamma(\mathscr X,\mathscr O_{\mathscr X}^\times) \to \mathscr C^\infty(\mathscr X(\C),\R) \to \hPic(\mathscr X)\to \Pic(\mathscr X) \to 0 ,
\end{equation}
where the first map is $f\mapsto \log \abs{f}$, 
the second associates with $\phi\in\mathscr C^\infty(\mathscr X(\C),\R)$
the trivial line bundle~$\mathscr O_{\mathscr X}$
endowed with the hermitian metric for which $\log \norm 1^{-1}=\phi$,
and the last one forgets the metric.

\subsection{}
The starting point of our lectures will be the following theorem
that asserts existence and uniqueness of ``arithmetic intersection degrees'' of
cycles associated with hermitian line bundles.
It fits naturally within
the  arithmetic intersection theory of~\cite{gillet-s90},
we refer to the foundational article by~\cite{bost-g-s94}
for such an approach; see also~\cite{faltings1992}
for a direct construction, as well as to the notes of Soulé in this volume.

\begin{theo}
Let $n=\dim(\mathscr X)$ and let 
$\overline{\mathscr L_1},\dots,\overline{\mathscr L_n}$ be
hermitian line bundles on~$\mathscr X$. There exists
a unique family of linear maps:
\[\hdeg  \left( \hc_1(\overline{\mathscr L_1})
   \cdots \hc_1(\overline{\mathscr L_d})\mid \cdot \right)
\colon \mathrm Z_d(\mathscr X)\to \R,
\]
for $d\in\{0,\dots,n\}$
satisfying the following  properties:
\begin{enumerate}
\item
For every integer~$d\in\{1,\dots,n\}$,
every integral closed subscheme~$Z$ of~$\mathscr X$ such that $\dim(Z)=d$,
every integer~$m\neq0$
and every regular meromorphic\footnote{that is,
defined over a dense open subscheme of~$Z$}
section~$s$ of~$\mathscr L_d^{m}|_Z$,
one has
\begin{multline}
m \hdeg \left( \hc_1(\overline{\mathscr L_1})
   \cdots \hc_1(\overline{\mathscr L_d})\mid Z \right)
\\ = \hdeg  \left( \hc_1(\overline{\mathscr L_1})
   \cdots \hc_1(\overline{\mathscr L_{d-1}})\mid \div(s) \right) \\
 {} + \int_{Z(\C)} \log\norm{s}^{-1} c_1(\overline{\mathscr L_1})
   \cdots c_1(\overline{\mathscr L_{d-1}}). 
\end{multline}
\item
For every closed point~$z$ of~$\mathscr X$,  viewed as a
integral closed subscheme of dimension~$d=0$, one has
\begin{equation}
 \hdeg\left( z \right)
 = \log (\Card( \kappa(z))),
\end{equation}
where $\kappa(z)$ is the residue field of~$z$.\footnote{It is indeed
a finite field.}
\end{enumerate}
Moreover, these maps are multilinear and symmetric in the hermitian
line bundles $\overline{\mathscr L_1},\dots,\overline{\mathscr L_n}$
and only depend on their isomorphism classes in~$\hPic(\mathscr X)$.
\end{theo}

\begin{rema}
This theorem should be put in correspondence with the analogous
geometric result for classical intersection numbers.
Let $F$ be a field and let $X$ be a proper scheme over~$F$,
let $n=\dim(X)$ and let $L_1,\dots,L_n$ be line bundles over~$X$.
The degree $ \deg ( c_1(L_1) \dots c_1(L_d)\mid Z) $
of a $d$-cycle~$Z$ in~$X$
is characterized by the relations:
\begin{enumerate}
\item It is linear in~$Z$;
\item If $d=0$ and $Z$ is a closed point~$z$ whose residue field $\kappa(Z)$
is a finite extension of~$F$, then 
$\deg (Z) = [\kappa(Z):F]$;
\item If $d\geq 1$ and $Z$ 
is an integral closed subscheme of~$X$ of dimension~$d$,
$m$ a non-zero integer, $s$ a regular meromorphic section of~$L_d^m|_Z$,
then 
\begin{equation}
 m \deg(c_1(L_1)\dots c_1(L_d)\mid Z)
 = \deg (c_1(L_1)\dots c_1(L_{d-1})\mid \div(s)). 
\end{equation}
\end{enumerate}
The additional integral that appears in the arithmetic degree
takes into account the fact that $\Spec(\Z)$  does not behave
as a proper variety.
\end{rema}

\begin{exem}
Assume that $Z$ is vertical and lies over a maximal ideal~$(p)$
of~$\Spec(\Z)$. Then $Z$ is a proper scheme over~$\F_p$ and
it follows from the inductive definition and the analogous
formula in classical intersection theory that
\[ 
 \hdeg \left( \hc_1(\overline{\mathscr L_1})
   \cdots \hc_1(\overline{\mathscr L_d})\mid Z \right)
 = \deg \left( c_1(\mathscr L_1|_{\mathscr X_{\F_p}})
   \cdots c_1({\mathscr L_d}|_{\mathscr X_{\F_p}})\mid Z \right) \log(p) .
\]
\end{exem}

\begin{exem}
Assume that $d=1$ and that $Z$ is horizontal, so that
$Z$ is the Zariski-closure in~$\mathscr X$ of a closed point
$z\in\mathscr X_\Q$. Let $F=\kappa(z)$ and let $\mathfrak o_F$
be its ring of integers; by properness of~$\mathscr X$,
the canonical morphism $\Spec(F)\to \mathscr X$ with image~$z$
extends to a morphism $\eps_z\colon\Spec(\mathfrak o_F)\to\mathscr X$,
whose image is~$Z$.
The formula
\[ \hdeg \left( \hc_1(\overline{\mathscr L})\mid Z\right)
 = \hdeg (\eps_z^* \overline{\mathscr L}) \]
expresses the arithmetic intersection number
as the \emph{arithmetic degree} of
the hermitian line bundle $\eps_z^*\overline{\mathscr L}$
over~$\Spec(\mathfrak o_F)$.
\end{exem}

\begin{prop}\label{prop.functoriality-degrees}
Let $f\colon \mathscr X'\to\mathscr X$ be a generically finite
morphism of proper flat schemes over~$\Z$,
let $Z$ be an integral closed subscheme of~$\mathscr X'$
and let $d=\dim(Z)$.

\begin{enumerate}
\item If $\dim(f(Z))<d$, then 
\[ \left( \hc_1(f^*\overline{\mathscr L_1})
   \cdots \hc_1(f^*\overline{\mathscr L_d})\mid Z \right) = 0; \]
\item Otherwise, $\dim(f(Z))=d$ and
\[ \left( \hc_1(f^*\overline{\mathscr L_1})
   \cdots \hc_1(f^*\overline{\mathscr L_d})\mid Z \right) = 
\left( \hc_1(\overline{\mathscr L_1})
   \cdots \hc_1(\overline{\mathscr L_d})\mid f_*(Z) \right),
 \]
where  $f_*(Z)=[ \kappa(Z):\kappa(f(Z))] f(Z)$
is a $d$-cycle on~$\mathscr X$.
\end{enumerate}
\end{prop}
This is the ``projection formula'' for height. It can be proved
by induction using the inductive definition of the heights
and the change of variables formula.

\begin{rema}
Let $n=\dim(\mathscr X)$ and assume that $\mathscr X$
is regular.
As the notation suggests rightly, the 
arithmetic intersection theory  of~\cite{gillet-s90}
allows another  definition of the real number
$ \hdeg \left( \hc_1(\overline{\mathscr L_1})
   \cdots \hc_1(\overline{\mathscr L_n})\mid \mathscr X \right)$
as the arithmetic degree of the $0$-dimensional arithmetic cycle
$ \hc_1(\overline{\mathscr L_1})
   \cdots \hc_1(\overline{\mathscr L_n}) \in \hCH_0(\mathscr X)$.

In fact, while the theory of~\cite{gillet-s90}
imposes regularity conditions on~$\mathscr X$, 
the definition of 
arithmetic product of classes of the form $\hc_1 (\overline{\mathscr L})$
requires  less stringent conditions;
in particular, the regularity of the
generic fiber $\mathscr X_\Q$ is enough.  
See~\cite{faltings1992} for such an approach, as well
as Soulé's notes in this volume.
More generally, for every birational morphism $f\colon Z'\to Z$
such that $Z'_{\Q}$ is regular, one
has
\[ \hdeg \left( \hc_1(\overline{\mathscr L_1})
   \cdots \hc_1(\overline{\mathscr L_d})\mid Z \right)
= 
\hdeg\left( \hc_1(f^*\overline{\mathscr L_1})
   \cdots \hc_1(f^*\overline{\mathscr L_d})\mid Z'\right). \]
\end{rema}

\section{The height of a variety}

\subsection{}
Let $X$ be a proper $\Q$-scheme and let $L$ be a line bundle on~$X$.
The important case is when the line bundle~$L$ is ample,
an assumption which will often be implicit below;
in that case, the pair $(X,L)$ is  called a polarized variety.

\subsection{}
Let $\mathscr X$ be a proper flat scheme over~$\Z$ and 
let $\overline{\mathscr L}$ be a hermitian line bundle on~$\mathscr X$
such that $\mathscr X_\Q=X$ and $\mathscr L_\Q=L$.
Let $Z$ be a closed integral subscheme of~$X$ and let $d=\dim(Z)$.
Let $\mathscr Z$ be the Zariski-closure of~$Z$ in~$\mathscr X$;
it is an integral closed subscheme of~$\mathscr X$ and $\dim(\mathscr Z)=d+1$.

\begin{defi}\label{defi.height}
The \emph{degree} and the \emph{height} of~$Z$ relative to~$\overline{\mathscr L}$
are defined by the formulas (provided $\deg_{\mathscr L}(Z)\neq0$).
\begin{gather}
\label{eq.defi-degree}\index{Degree>of a subscheme}\glossary{$\deg_{\mathscr L} (Z)$: Degree}
\deg_{\mathscr L} (Z) = \deg ( c_1(L)^d \mid Z ) \\
\label{eq.defi-height}\index{Height>of a subscheme}\glossary{$h_{\overline{\mathscr L}}(Z)$: Height}
h_{\overline{\mathscr L}}(Z) = \hdeg\left( \hc_1(\overline{\mathscr L}^{d+1})\mid \mathscr Z\right) / \big( (d+1) \deg_{\mathscr L}(Z) \big) .
\end{gather}
\end{defi}

Note that the degree $\deg_{\mathscr L}(Z)$ is computed on~$X$, 
hence only depends  on~$L$.  Moreover, 
the condition that $\deg_{\mathscr L}(Z)\neq 0$ is satisfied
(for every~$Z$) when $L$ is ample on~$X$.

\begin{prop}\label{prop.functoriality-heights}
Let $f\colon\mathscr X'\to\mathscr X$ be a generically finite morphism
of proper flat schemes over~$\Z$, 
let $Z$ be a closed integral subscheme of~$\mathscr X'_\Q$
and let $d=\dim(Z)$. Assume that $L$ is ample on~$X$ and that
$\dim(f(Z))=d$. Then $\deg_{f^*\mathscr L}(Z)>0$ and
\[ h_{f^*\overline{\mathscr L}}(Z) = h_{\overline{\mathscr L}}(f(Z)). \]
\end{prop}
\begin{proof}
This follows readily from proposition~\ref{prop.functoriality-degrees}
and its analogue for geometric degrees. Indeed, when one compares
 formula~\eqref{eq.defi-height} for $Z$ and for~$f(Z)$, both the numerator 
and the denominator
get multiplied by $[\kappa(Z):\kappa(f(Z))]$.
\end{proof}

\begin{exem}
For every $x\in X(\overline\Q)$, let $[x]$ denote its Zariski closure in~$X$.
The function $X(\overline\Q)\to \R$ given by $x\mapsto h_{\overline{\mathscr L}}([x])$ is a height function relative to the line bundle~$\mathscr L_\Q$ on~$X$
in the sense of Weil. Note that while the height functions
defined by  Weil which are associated with a line bundle~$L$ on~$X$
are only defined  modulo the space of bounded  functions on~$X(\overline\Q)$,
the techniques of Arakelov geometry show that specifying
a model~$\mathscr X$ of~$X$ and a hermitian line bundle~$\overline{\mathscr L}$
on~$\mathscr X$ such that $\mathscr L_\Q$ is isomorphic to~$L$
allows to obtain a well-defined representative of this class,
for which precise theorems, such as the equidistribution theorem below,
can be proved.
\end{exem}

\begin{exem}\label{exem.neron-tate-good}
Let us assume that $X$ is an abelian variety over a number field~$F$,
with everywhere good reduction,
and let $\mathscr X$ be an $\mathfrak o_F$-abelian scheme such 
that $\mathscr X_F=X$.
Let $o$ be the origin of~$X$ and let $\eps_o\colon\Spec(\mathfrak o_F)\to\mathscr X$ be the corresponding section.
Let $L$ be a line bundle on~$X$ with a trivialisation~$\ell$ of~$L|_o$.
There exists a unique line bundle~$\mathscr L$ on~$\mathscr X$
such that $\mathscr L_F=L$ and
such that the given trivialisation of~$L|_o$ extends to a trivialisation
of $\eps_o^*\mathscr L$.
Moreover, for every embedding~$\sigma\colon F\hra\C$
the theory of Riemann forms on complex tori endows~$L_\sigma$
with a canonical metric $\norm\cdot_\sigma$ 
whose curvature form $c_1(L_\sigma,\norm\cdot_\sigma)$
is invariant by translation
and such that $\norm{\ell}_\sigma=1$; this is in fact the unique metric
possessing these two properties.
We let $\overline{\mathscr L}$ be the hermitian line bundle on~$\mathscr X$
so defined.  

The associated height function will be denoted by~$\widehat h_L$:
it extends the \emph{Néron--Tate height}\index{Nerontateheight=N\'eron--Tate height}\index{Height>Nerontate=(N\'eron--Tate ---)}\glossary{$\widehat h_L$: N\'eron-Tate height} from~$X(\overline\Q)$ 
to all integral closed subschemes.

Assume that $L$ is even, that is $[-1]^*L\simeq L$.
Then $[n]^*L \simeq L^{n^2}$ for every integer~$n\geq 1$, and this
isomorphism extends to an isomorphism of hermitian line bundles
$[n]^*\overline{\mathscr L}\simeq \overline{\mathscr L}^{n^2}$.
Consequently, for every integral closed subscheme~$Z$ of~$X$, one has
the following relation
\begin{equation}
\widehat h_L([n](Z)) = n^2 \widehat h_L(Z). 
\end{equation}

Assume otherwise that $L$ is odd, that is $[-1]^*L\simeq L^{-1}$.
Then $[n]^*L\simeq L^{n}$ for every integer~$n\geq 1$; similarly, this
isomorphism extends to an isomorphism of hermitian line bundles
$[n]^*\overline{\mathscr L}\simeq \overline{\mathscr L}^{n}$.
Consequently, for every integral closed subscheme~$Z$ of~$X$, one has
the following relation
\begin{equation}
\widehat h_L([n](Z)) = n \widehat h_L(Z). 
\end{equation}
\end{exem}

\begin{prop}
Let $\mathscr X'$ be a proper flat scheme over~$\Z$
such that $\mathscr X'_\Q=X$; let $\overline{\mathscr L'}$
be a hermitian line bundle on~$\mathscr X'$ such that $\mathscr L'_\Q=L$.
Assume that $L$ is ample. 
Then there exists a real number~$c$ such that
\[ \abs{h_{\overline {\mathscr L}}(Z) - h_{\overline{\mathscr L'}}(Z)}
\leq c \]
for every integral closed subscheme~$Z$ of~$X$.
\end{prop}
\begin{proof}
One proves in fact the existence of a real number~$c$ such that
\[
\abs{ \hdeg(\hc_1(\overline{\mathscr L})^{d+1}\mid \mathscr Z)
-  \hdeg(\hc_1(\overline{\mathscr L'})^{d+1}\mid \mathscr Z')}
\leq
 c \deg( c_1(L)^{d}\mid Z) \]
for every integral $d$-dimensional subvariety~$Z$ of~$X$,
where $\mathscr Z$ and $\mathscr Z'$ are the Zariski closures
of~$Z$ in~$\mathscr X$ and~$\mathscr X'$ respectively.
Considering a model $\mathscr X''$ that dominates
$\mathscr X$ and $\mathscr X'$ (for example, the Zariski closure
in $\mathscr X\times_\Z \mathscr X'$ of the diagonal),
we may assume that $\mathscr X=\mathscr X'$, hence
$\mathscr Z=\mathscr Z'$.
A further reduction, that we omit here, allows us to assume
that $\mathscr L$ is a nef line bundle, and that its hermitian
metric is semipositive, and similarly for~$\mathscr L'$.

By multilinearity, the left hand side that we wish to bound from
above is the absolute value of 
\[ \sum_{i=0}^{d} \hdeg( \hc_1(\overline{\mathscr L}'\otimes\overline {\mathscr L}^{-1}) \hc_1(\overline{\mathscr L}')^i \hc_1(\overline{\mathscr L})^{d-i}\mid \mathscr Z). \]
Then we view the section~$1$ of $\mathscr L'_\Q\otimes \mathscr L^{-1}_\Q$
as a meromorphic section~$s$ of $\mathscr L'\otimes\mathscr L^{-1}$.
Note that its divisor is purely vertical, and its hermitian
norm $\norm s$ is a non-vanishing continuous function on~$X(\C)$.
The definition of the arithmetic intersection numbers
then leads us to estimate algebraic intersection numbers
\[  \deg( c_1(\mathscr L')^i c_1(\mathscr L)^{d-i}
\mid \div(s|_{ \mathscr Z})) \]
and an integral
\[ \int_{Z(\C)}  \log(\norm{s}^{-1}) 
c_1(\overline{\mathscr L}')^i c_1(\overline{\mathscr L})^{d-i}. \]

By positivity of the curvatures
forms $c_1(\overline{\mathscr L})$ and $c_1(\overline{\mathscr L'})$,
the latter integral is bounded from above by
\[
\left\lVert {\log(\norm{s})^{-1}} \right\rVert_\infty
\int_{Z(\C)} c_1(\overline{\mathscr L}')^i c_1(\overline{\mathscr L})^{d-i}
=  \left\lVert {\log(\norm{s})^{-1}} \right\rVert_\infty
 \deg( c_1(L)^d\mid Z).
\]

The algebraic terms can be bounded as well.
Observe that there exists an integer~$n\geq 1$ such that~$ns$ extends
to a global section of $\mathscr L'\otimes\mathscr L^{-1}$,
and $n s^{-1}$ extends
to a global section of its inverse $(\mathscr L')^{-1}\otimes\mathscr L)$.
(This is the ultrametric counterpart to the fact that the section~$s$
has non-vanishing norm on~$X(\C)$.)
Consequently, $\div(ns|_{\mathscr Z})$ and $\div(ns^{-1}|_{\mathscr Z})$
are both effective, so that 
\[ - \sum_p v_p(n) [\mathscr Z_{\F_p}] \leq \div(s|_{\mathscr Z}) \leq
 \sum_p v_p(n) [\mathscr Z_{\F_p}]. \]
This inequality of cycles is preserved after taking
intersections, so that
\begin{align*} 
\deg(c_1(\mathscr L')^ic_1(\mathscr L)^{d-i}\mid \div(s|_{\mathscr Z})_p )
\hskip -3cm \\
& \leq
v_p(n) \deg(c_1(\mathscr L')^ic_1(\mathscr L)^{d-i}\mid [\mathscr Z_{\F_p}])\\
& = 
v_p(n) \deg(c_1(L)^d\mid Z),\end{align*}
where $\div(s|_{\mathscr Z})_p$ 
is the part of~$\div(s|_{\mathscr Z})$ that lies above the maximal
ideal~$(p)$ of~$\Spec(\Z)$.
There is a similar lower bound.

Adding all these contributions, this proves the proposition.
We refer to~\cite{bost-g-s94}, \S3.2.2, for more details.
\end{proof}


\begin{prop}
Let us assume that~$L$ is ample.
For every real number~$B$,
the set of integral closed subschemes~$Z$ of~$X$
such that $\deg_L(Z)\leq B$ and $h_{\overline{\mathscr L}}(Z)\leq B$
is finite.
\end{prop}
The case of closed points is Northcott's theorem,
and the general case is Theorem~3.2.5 of~\cite{bost-g-s94}.
The principle of its proof goes by reducing
to the case where $X=\P^N$ and $\overline{\mathscr L}=\overline{\mathscr O(1)}$,
and comparing the height~$h_{\overline{\mathscr L}}(Z)$ of
a closed integral subscheme~$Z$
with the height of its Chow form.
(That paper also provides a more elementary proof,
relying on the fact that a finite set of sections of powers
of~$\mathscr O(1)$ are sufficient to compute
by induction the height of any closed integral subscheme
of~$\mathbf P^N$ of given degree.)

\section{Adelic metrics}

\subsection{}
Let $S=\{2,3,\dots,\infty\}$ be the set of places of~$\Q$.

Each prime number~$p$ is identified with the $p$-adic absolute value on~$\Q$,
normalized by $\abs{p}_p=1/p$; these places are said to be finite.
We denote by $\Q_p$ the completion of~$\Q$ for this $p$-adic absolute value
and fix an algebraic closure $\overline{\Q_p}$ of~$\Q_p$.
The $p$-adic absolute value extends uniquely to~$\overline{\Q_p}$;
the corresponding completion is denoted by~$\C_p$: this is an algebraically closed  complete valued field.

The archimedean place is represented
by the symbol~$\infty$, and is identified with
the usual absolute value on~$\Q$; it is also called the infinite place.
For symmetry of notation, we may write $\Q_\infty=\R$ and $\C_\infty=\C$,
the usual fields of real and complex numbers.

The notion of adelic metrics that we now describe is close in spirit to
the one exposed by Peyre in this volume. However, Peyre is only
interested in rational points and it suffices for him to consider
local fields such as~$\Q_p$, and ``naïve $p$-adic analytic varieties''.
We will however need to consider
algebraic points. The topological spaces associated with varieties
over~$\C_p$ are not well behaved enough, and we will have to
use a rigid analytic notion of $p$-adic varieties, specifically
the one introduced by Berkovich.

\subsection{}\label{ss.berkovich}
Let $X$ be a proper scheme over~$\Q$.
Let $v\in S$ be a place of~$\Q$.

Assume $v=\infty$. Then we set $X_\infty^\an=X(\C_\infty)/F_\infty$, the set of complex points of~$X$ modulo the action of complex conjugation~$F_\infty$.

Assume now that $v=p$ is a finite place. Then we set $X_p^\an$
to be the analytic space  associated by~\cite{berkovich1990}\index{Berkovich space}\index{Space>berkovich=(Berkovich ---)}\glossary{$X_p^{\mathrm{an}}$: Berkovich space}
to the $\Q_p$-scheme $X_p=X_{\Q_p}$. 
As a set, $X_p^\an$ is the quotient of the class of all pairs
$(K,x)$, where $K$ is a complete valued field containing~$\Q_p$
and $x\in X(K)$, modulo the identification $(K,x)\sim (K',x')$
if there exists a common valued field extension~$K''$ of~$K$ and~$K'$
such that $x$ and~$x'$ induce the same point of~$X(K'')$.
Its topology is the coarsest topology such that  for every open 
subscheme~$U$ of~$X$ and every $f\in\mathscr O_X(U)$,
the subset $U_p^\an$ (consisting of classes $[K,x]$ such that $x\in U(K)$)
is open in~$X_p^\an$,
and the function $[K,x]\mapsto \abs{f(x)}$ from~$U_p^\an$ to~$\R_+$
is continuous.

The space~$X_p^\an$ is a compact metrizable topological space, 
locally contractible (in particular locally arcwise connected).
There is a canonical continuous map $X(\C_p)\to X_p^\an$;
it identifies the (totally discontinuous) topological space 
$X(\C_p)/\Gal(\C_p/\Q_p)$ with a dense subset of $X_p^\an$.
It is endowed with a sheaf in local rings~$\mathscr O_{X_p^\an}$;
for every open subset~$U$ of~$X_p^\an$, every holomorphic function
$f\in\mathscr O_{X_p^\an}(U)$ admits an absolute value
$\abs f\colon U\to\R_+$.

We gather all places together and consider
the topological space $X_\ad=\coprod_{v\in S}X_v^\an$,
coproduct of the family $(X_v^\an)_{v\in S}$.
By construction, a function $\phi$ on~$X_\ad$ consists in a family
$(\phi_v)_{v\in S}$, where $\phi_v$ is a function on~$X_v^\an$,
for every $v\in S$.

\subsection{}
Let $L$ be a line bundle on~$X$; it induces a line bundle~$L_v^\an$
on~$X_v^\an$ for every place~$v$.

A continuous $v$-adic metric on~$L_v^\an$ is the datum,
for every open subset~$U$ of~$X_v^\an$ and every section~$s$ of~$L_v^\an$
on~$U$, of a continuous function $\norm s\colon U\to \R_+$,
subject to the requirements: 
\begin{enumerate}
\item For every open subset~$V$ of~$U$, one has $\norm {s|_V}=\norm{s}|_V$;
\item For every holomorphic function $f\in\mathscr O_{X_v^\an}(U)$,
one has $\norm{fs}=\abs f \norm s$.
\item If $s$ does not vanish, then $\norm s$ does not vanish as well.
\end{enumerate}

If $L$ and~$M$  are line bundles on~$X$ equipped
with $v$-adic metrics, then $L^{-1}$ and $L\otimes M$ admit
natural $v$-adic metrics, and the canonical isomorphism 
$L^{-1}\otimes L\simeq\mathscr O_X$ is an isometry.

The trivial line bundle $\mathscr O_X$ admits a canonical $v$-adic metric
for which $\norm f=\abs f$ for every local section of~$\mathscr O_X$.
More generally, for every $v$-adic metric $\norm\cdot$ on~$\mathscr O_X$,
$\phi=\log \norm 1^{-1}$ is a continuous function on~$X_v^\an$,
and any $v$-adic metric on~$\mathscr O_X$ is of this form.
The $v$-adically metrized line bundle associated with~$\phi$
is denoted by~$\mathscr O_X(\phi)$.

If $\overline L$ is a line bundle endowed with an $v$-adic metric
and $\phi\in \mathscr C(X_v^\an,\R)$, we denote by $\overline L(\phi)$
the $v$-adically metrized line bundle $\overline L \otimes \mathscr O_X(\phi)$.
Explicitly, its $v$-adic metric is that of~$\overline L$ multiplied
by~$e^{-\phi}$.

\begin{exem}
Let $\mathscr X$ be a proper flat scheme over~$\Z$ such that $\mathscr X_\Q=X$,
let $d$ be a positive integer
and let ${\mathscr L}$ be a line bundle on~$\mathscr X$
such that $\mathscr L_\Q=L^d$. Let us show that this datum endows~$L$
with an $p$-adic metric, for every finite place~$p\in S$.


Let thus fix a prime number~$p$.
There exists a canonical \emph{specialization map},
$X_p^\an\to \mathscr X\otimes_\Z\F_p$; it is anticontinuous
(the inverse image of an open subset is closed).
For every open subset~$\mathscr U\subset\mathscr X\otimes_\Z\F_p$, 
let $\tube{\mathscr U}$ be the preimage of~$\mathscr U$.

There exists a unique continuous metric on~$L_p^\an$ such that
for every open subscheme $\mathscr U$ of~$\mathscr X\otimes_\Z\Z_p$
and every basis $\ell$ of~$\mathscr L$ on~$\mathscr U$,
one has $\norm\ell\equiv 1$ on $\tube{\mathscr U\otimes\F_p}$.
Explicitly, if $s$ is a section of~$L_p^\an$ on an open subset~$U$
of~$\tube{\mathscr U\otimes\F_p}$,
there exists a holomorphic function $f\in\mathscr O_{X_p^\an}(U)$
such that $s^d=f \ell$ and $\norm s =\abs{f}^{1/d}$ on~$U$.

Such $p$-adic metrics are called \emph{algebraic}.\index{Algebraic metric}
\end{exem}

\subsection{}
An adelic metric on~$L$ is the datum, for every place $v\in S$,
of a $v$-adic metric on the line bundle~$L_v^\an$ on~$X_v^\an$,
subject to the additional requirement that there exists
a model $(\mathscr X,\mathscr L)$ of $(X,L)$ inducing
the given $p$-adic metric for all but finitely many primes~$p$.

If $L$ and~$M$  are line bundles on~$X$ equipped
with adelic metrics, then $L^{-1}$ and $L\otimes M$ admit
natural adelic metrics, and the canonical isomorphism 
$L^{-1}\otimes L\simeq\mathscr O_X$ is an isometry.

The trivial line bundle $\mathscr O_X$ admits a canonical adelic metric
for which $\norm f=\abs f$ for every local section of~$\mathscr O_X$.
More generally, for every adelic metric $\norm\cdot$ on~$\mathscr O_X$,
and every place~$v\in S$, 
then $\phi_v=\log \norm 1_v^{-1}$ is a continuous function on~$X_v^\an$,
and is identically zero for all but finitely many places~$v$;
in other words, the function $\phi=(\phi_v)\in\mathscr C(X_\ad,\R)$
has compact support.
The adelically metrized line bundle associated with~$\phi$
is denoted by~$\mathscr O_X(\phi)$.
Conversely, any adelic metric on~$\mathscr O_X$ is of this form.

If $\overline L$ is a line bundle endowed with an adelic metric
and $\phi\in \mathscr C_\cpct(X_\ad,\R)$, we denote by $\overline L(\phi)$
the adelically metrized line bundle $\overline L \otimes \mathscr O_X(\phi)$.
Explicitly, for every place~$v$,
its $v$-adic metric is that of~$\overline L$ multiplied by~$e^{-\phi_v}$.

\begin{rema}\label{rema.all-but-finitely}
Let $(\mathscr X,\mathscr L)$ and $(\mathscr X',\mathscr L')$
be two models of the polarized variety $(X,L)$.
Since $X$ is finitely presented, there exists a dense open subscheme~$U$
of~$\Spec(\Z)$ such that the isomorphism $\mathscr X_\Q=X=\mathscr X'_\Q$
extends to an isomorphism $\mathscr X_U\simeq \mathscr X'_U$.
Then, up to shrinking~$U$,
we may assume that the isomorphism $\mathscr L_\Q=L=\mathscr L'_\Q$
extends to an isomorphism $\mathscr L_U\simeq \mathscr L'_U$.
In particular, for every prime number~$p$ such that $(p)\in U$,
the $p$-adic norms on~$L$ induced by~$\mathscr L$ and~$\mathscr L'$
coincide.
\end{rema}

\subsection{}
Let $\bPic(X_\ad)$ be the abelian group of isomorphism
classes of line bundles on~$X$ endowed with adelic metrics. It fits within an exact sequence
\begin{equation}\label{eq.pic}
  \Gamma(X,\mathscr O_X^\times) \to \mathscr C_\cpct(X_\ad,\R) \to \bPic(X_\ad) \to \Pic(X) \to 0. \end{equation}
The morphism on the left is given by $u\mapsto (\log\abs u_v^{-1})_{v\in S}$.
It is injective up to torsion, as a consequence of Kronecker's theorem:
if $\abs u_v=1$ for every place~$v$, 
then there exists~$m\geq 1$ such that $u^m= 1$.
Its image is the kernel of the morphism $\mathscr C(X_\ad,\R)\to \bPic(X)$; indeed,
an isometry $\mathscr O_X(\phi)\to \mathscr O_X(\psi)$ is given
by an element $u\in \Gamma(X,\mathscr O_X^\times)$
such that $\psi_v+\log\abs u_v^{-1}=\phi_v$, for every place $v\in S$.

We denote by $\hc_1(\overline L)$ the isomorphism class
in $\bPic(X_\ad)$ of an adelically metrized line bundle~$\overline L$ on~$X$.

\begin{rema}
Let $D$ be an effective Cartier divisor on~$X$ and let $\mathscr O_X(D)$
be the corresponding line bundle; let $s_D$ be its canonical section.
Assume that $\mathscr O_X(D)$ is endowed with an adelic metric.

Let $v\in S$ be a place of~$\Q$.
The function $g_{D,v}=\log\norm{s_D}_v^{-1}$
is a continuous function on $X_v^\an\setminus \abs D$, and is called
a $v$-adic Green function for~$D$.
For every open subscheme~$U$ of~$X$ and any equation~$f$
of~$D$ on~$U$, $g_{D,v}+\log \abs f_v$ extends to a continuous function
on $U_v^\an$. 
Conversely, this property characterizes $v$-adic Green functions for~$D$.

The family $g_D=(g_{D,v})$ is called an adelic Green function for~$D$.
\end{rema}

\begin{lemm}[\cite{chambert-loir-thuillier2009}, prop.~2.2]
Let $\mathscr X$ be a proper flat integral scheme over~$\Z$, 
let $\overline{\mathscr L}$
be a hermitian line bundle on~$\mathscr X$. Let $X=\mathscr X_\Q$
and let $L=\mathscr L_\Q$, endowed with the algebraic adelic metric
associated with $(\mathscr X,\overline{\mathscr L})$.
Assume that $\mathscr X$ is integrally closed  
in its generic fiber (for example, that it is normal).

Then the canonical map
$\Gamma(\mathscr X,\mathscr L)\to \Gamma(X,L)$
is injective and its image is the set of sections~$s$
such that $\norm s_v\leq 1$  for every finite place~$v\in S$.
\end{lemm}
Equivalently, effective Cartier divisors on~$\mathscr X$ correspond
to $v$-adic Green functions which are positive at all finite places~$v$.
\begin{proof}
Injectivity follows from the fact that $\mathscr X$ is flat,
so that $X$ is schematically dense in~$\mathscr X$.
Surjectivity is a generalization of the fact that an integrally closed
domain is the intersection of its prime ideals of height~1.
\end{proof}

\subsection{}
Let $\norm\cdot$ and $\norm\cdot'$ be two adelic metrics on~$L$.
The ratio of these metrics is a metric on the trivial line bundle,
hence is of the form $\mathscr O_X(\phi)$,
for some function $\phi=(\phi_v)\in\mathscr C_\cpct(X_\ad,\R)$.
For all but finitely many places, $\phi_v$ is identically~$0$.
For every place~$v$, we let 
\[ \delta_v(\norm\cdot,\norm\cdot')= \norm{\phi_v}_\infty = \sup_{x\in X_v^\an} \left| \log \frac{\norm{\cdot}'}{\norm\cdot} (x) \right|. \]
Since $X_v^\an$ is compact, this is a positive real number;
it is equal to~$0$ for all but finitely many places~$v$.

We then define  the distance between the two given adelic metrics by
\[ \delta (\norm\cdot,\norm\cdot')= \sum_{v\in S} \delta_v (\norm{\cdot}',\norm\cdot).\]

The set of adelic metrics on a given line bundle~$L$ is a real affine space,
its underlying vector space  is the subspace $\mathscr C_\cpct(X_\ad,\R)$
of $\mathscr C(X_\ad,\R)=\prod_v \mathscr C(X_v^\an,\R)$ 
consisting of families $(\phi_v)$
such that $\phi_v\equiv 0 $ for all but finitely many places~$v\in S$.

The space $\mathscr C_\cpct(X_\ad,\R)$  is the union
of the subspaces $\mathscr C_U(X_\ad,\R)$
of functions with (compact) support above a given finite set~$U$ of places of~$S$.
We thus endow it with its natural inductive limit topology.

\begin{exem}[Algebraic dynamics]\label{exem.algebraic-dynamics}
We expose here the point of view of~\cite{zhang95b}
regarding the canonical heights of~\cite{CallSilverman-1993}
in the context of algebraic dynamics.
Let $X$ be a proper $\Q$-scheme, 
let $f\colon X\to X$ be a morphism,
let $L$ be a line bundle on~$X$ such that $f^*L\simeq L^q$,
for some integer~$q\geq 2$. We fix such an isomorphism~$\eps$.
The claim is that \emph{there exists a unique adelic metric on~$L$
for which the isomorphism~$\eps$ is an isometry.}

Let us first fix a place~$v$ and prove that there is a unique
$v$-adic metric on~$L$  for which $\eps$ is an isometry.
To that aim, let us consider, for any $v$-adic metric $\norm\cdot$ on~$L$,
the induced $v$-adic metric on~$f^*L$
and transfer it to~$L^q$ via~$\eps$. This furnishes a $v$-adic
metric $\norm\cdot^f$ on~$L$ such that $\eps$ is an isometry
from $(L,f^*\norm \cdot)$ to $(L, \norm\cdot ^f)^q$, and it 
is the unique $v$-adic metric on~$L$ satisfying this property.
Within the real affine space of $v$-adic metrics on~$L$,
normed by the distance~$\delta_v$, and complete,
the self-map $\norm\cdot \mapsto\norm\cdot^f$ is contracting
with Lipschitz constant~$1/q$.  Consequently, the claim
follows from Banach's fixed point theorem.

We also note that there exists a dense open subscheme~$U$
of $\Spec(\Z)$, a model $(\mathscr X,\mathscr L)$
of $(X,L)$ over~$U$  such that $f\colon X\to X$ extends to
a morphism $\phi\colon\mathscr X\to\mathscr X$
and the isomorphism $\eps\colon f^*L\simeq L^q$ extends to an isomorphism
$\phi^*\mathscr L\simeq \mathscr L^q$, still denoted by~$\eps$.
This implies that for every finite place~$p$ above~$U$,
the canonical $v$-adic metric is induced 
by the model~$(\mathscr X,\mathscr L)$.

Consequently, the family  $(\norm\cdot_v)$ of $v$-adic metrics on~$L$
for which $\eps$ is an isometry is an adelic metric.
\end{exem}

\section{Arithmetic ampleness}


\begin{defi}
Let $\mathscr X$ be a proper scheme over~$\Z$ and let $\overline{\mathscr L}$
be a hermitian line bundle on~$\mathscr X$.
One says that $\overline{\mathscr L}$ is relatively \emph{semipositive}\index{Semipositive>hermitian line bundle}
if:
\begin{enumerate}
\item For every vertical integral curve~$C$ on~$\mathscr X$, one has
$\deg_{\mathscr L}(C)\geq0$;
\item For every holomorphic map $f\colon \mathbf D\to \mathscr X(\C)$,
the curvature of $f^*\overline{\mathscr L}$ is semipositive.
\end{enumerate}
\end{defi}

If $\overline{\mathscr L}$ is relatively semipositive,
then $\mathscr L_\Q$ is nef. (The degree of a curve~$Z$ in~$X_\Q$
with respect to~$L_\Q$ of the integral over~$Z(\C)$ of the curvature 
of $\overline{\mathscr L}$, hence is positive.)

\begin{exem}
Let us consider the tautological line bundle $\mathscr O(1)$
on the projective space $\P^N_\Z$. Its local sections
correspond to homogeneous rational functions of
degree~$1$ in indeterminates $T_0,\dots,T_N$. 
It is endowed with a natural hermitian metric
such that, if $f$ is such 
a rational function, giving rise to the section~$s_f$,
and if $x=[x_0:\dots:x_N]\in\P^N(\C)$, one has the formula
\[ \norm{s_f}(x) = \frac{\abs{f(x_0,\dots,x_N)}}{\left(\abs{x_0}^2+\dots+\abs{x_N}^2\right)^{1/2}}. \]
(By homogeneity of~$f$, the right hand side does not depend
on the choice of the system of homogeneous coordinates for~$x$.)
The corresponding hermitian line bundle $\overline{\mathscr O(1)}$
is relatively semipositive.
It is in fact the main source of relatively semipositive hermitian line bundles,
in the following way.

Let $\mathscr X$ be a proper scheme over~$\Z$ and let $\overline{\mathscr L}$
be a hermitian line bundle on~$\mathscr X$.
One says that $\overline{\mathscr L}$ is relatively \emph{ample}\index{Ample>hermitian line bundle}
if there exists an embedding $\phi\colon\mathscr X\hra \P^N_\Z$,
a metric with strictly positive curvature on~$\mathscr O_{\P^N}(1)$
and an integer~$d\geq1$
such that $\overline{\mathscr L}^d\simeq \phi^* \overline{\mathscr O_{\P^N}(1)}$.
\end{exem}

\begin{prop}\label{prop.continuity-height}
Let $X$ be a proper scheme over~$\Q$ and 
let $L_0,\dots,L_d$ be line bundles on~$X$.

Let $\mathscr X,\mathscr X'$ be proper flat schemes over~$\Z$ such that
$X=\mathscr X_\Q=\mathscr X'_\Q$, 
let $\overline{\mathscr L_0},\dots,\overline{\mathscr L_d}$ 
(respectively $\overline{\mathscr L'_0},\dots,\overline{\mathscr L'_d}$)
be 
semipositive hermitian line bundles on~$\mathscr X$ (respectively $\mathscr X'$)
such that $\mathscr L_{j,\Q}=\mathscr L'_{j,\Q}=L_j$.
We write $\overline{L_0},\dots,\overline{L_d}$ (respectively $\overline{L'_0},\dots,\overline{L'_d}$)
for the corresponding adelically metrized line bundles on~$X$. 

Then for every closed subscheme~$Z$ of~$\mathscr X$, one has
\begin{multline*}
 \Big| \hdeg\left( \hc_1(\overline{\mathscr L'_0})\dots\hc_1(\overline{\mathscr L'_d})\mid Z \right)
-  \hdeg\left( \hc_1(\overline{\mathscr L_0})\dots\hc_1(\overline{\mathscr L_d})\mid Z \right) \Big| \\
\leq \sum_{j=0}^d \delta(\overline{L_j},\overline{L'_j})   \deg \left( c_1(L_0)\dots \widehat{c_1(L_j)} \dots c_1(L_d)\mid Z \right),
\end{multline*}
where the factor $c_1(L_j)$ is omitted in the $j$th term.
\end{prop}
\begin{proof}
Using the projection formula, 
we first reduce to the case where $\mathscr X=\mathscr X'$ is normal. 
We then write
{\small
\begin{multline*}
\hdeg\left( \hc_1(\overline{\mathscr L'_0})\dots\hc_1(\overline{\mathscr L'_d})\mid Z \right) 
-  \hdeg\left( \hc_1(\overline{\mathscr L_0})\dots\hc_1(\overline{\mathscr L_d})\mid Z \right)  \\
= \sum_{j=0}^d
\hdeg\left( \hc_1(\overline{\mathscr L'_0})\dotsc \hc_1(\overline{\mathscr L'_{j-1}})(\hc_1(\overline{\mathscr L'_j})- \hc_1(\overline{\mathscr L_j}))
\hc_1(\overline{\mathscr L_{j+1}})\dotsc\hc_1(\overline{\mathscr L_d})\mid Z \right) 
\end{multline*}}
and bound the $j$th term as follows.
Let $s_j$ be the regular meromorphic section
of $\mathscr O_X=\mathscr L'_j\otimes (\mathscr L_j)^{-1}$ corresponding to~$1$.
By definition, one has
{\small
\begin{align*}
\hdeg\left( \hc_1(\overline{\mathscr L'_0})\dots \hc_1(\overline{\mathscr L'_{j-1}})(\hc_1(\overline{\mathscr L'_j})- \hc_1(\overline{\mathscr L_j}))
\hc_1(\overline{\mathscr L_{j+1}})\dots\hc_1(\overline{\mathscr L_d})\mid Z \right) \hskip -10cm \\
& =\hdeg\left( \hc_1(\overline{\mathscr L'_0})\dots \hc_1(\overline{\mathscr L'_{j-1}})
\hc_1(\overline{\mathscr L_{j+1}})\dots\hc_1(\overline{\mathscr L_d})
(\hc_1(\overline{\mathscr L'_j})- \hc_1(\overline{\mathscr L_j}))
\mid Z \right)  \\
& =\hdeg\left( \hc_1(\overline{\mathscr L'_0})\dots \hc_1(\overline{\mathscr L'_{j-1}})
\hc_1(\overline{\mathscr L_{j+1}})\dots\hc_1(\overline{\mathscr L_d})
\mid \div(s_j|_Z) \right)  \\
& \qquad + \int_{Z(\C)} \log \norm{s_j}^{-1} c_1(\overline{\mathscr L'_0})\dots c_1(\overline{\mathscr L'_{j-1}})
c_1(\overline{\mathscr L_{j+1}})\dots c_1(\overline{\mathscr L_d})
.
\end{align*}}
Moreover, all components  of~$\div(s_j|_Z)$  are vertical.
For every $j\in\{0,\dots,d\}$ and every $v\in S$,
let $\delta_{j,v}= \delta_v(\overline {L_j},\overline {L'_j})$
(this is zero for all but finitely many places~$v$).
Using the fact that
$\abs{\log\norm{s_j}_v } \leq \delta_v(\overline {L_j},\overline {L'_j})$
for every place~$v\in S$,
the normality assumption on~$\mathscr X$ 
implies that  
\[
\div(s_j|_Z) \leq  \sum_{p\in S\setminus\{\infty\}}\delta_{j,p} (\log p)^{-1} [Z\otimes\F_p].
\]
Since the line bundles~$\mathscr L_k$ and $\mathscr L'_k$ are semipositive,
this implies the bound
{\small
\begin{align*}
\hdeg\left( \hc_1(\overline{\mathscr L'_0})\dots \hc_1(\overline{\mathscr L'_{j-1}})
\hc_1(\overline{\mathscr L_{j+1}})\dots\hc_1(\overline{\mathscr L_d})
\mid \div(s_j|_Z) \right)  \hskip -8cm \\
& =
\sum_{p} \deg\left( c_1({\mathscr L'_0})\dots c_1({\mathscr L'_{j-1}})
c_1({\mathscr L_{j+1}})\dots c_1({\mathscr L_d})
\mid \div(s_j|_Z)_p \right) \log p  \\
& \leq \sum_p \delta_{j,p} 
\deg\left( c_1({\mathscr L'_0})\dots c_1({\mathscr L'_{j-1}})
c_1({\mathscr L_{j+1}})\dots c_1({\mathscr L_d})
\mid [Z\otimes\F_p]  \right)   \\
& \leq \left( \sum_p \delta_{j,p}  \right)
\deg\left( c_1(L_0)\dots c_1(L_{j-1})
c_1(L_{j+1})\dots c_1(L_d) \mid Z \right). 
\end{align*}}
Similarly, the curvature forms $c_1(\overline{L_k})$ and $c_1(\overline{L'_k})$ are semipositive, so that the upper bound
 $\log\norm{s_j}^{-1}\leq \delta_{j,\infty}$
implies
\begin{align*}
\int_{Z(\C)} \log \norm{s_j}^{-1} c_1(\overline{\mathscr L'_0})\dots c_1(\overline{\mathscr L'_{j-1}})
c_1(\overline{\mathscr L_{j+1}})\dots c_1(\overline{\mathscr L_d}) 
 \hskip-7cm \\
& \leq 
\delta_{j,\infty}
 \int_{Z(\C)}  c_1(\overline{\mathscr L'_0})\dots c_1(\overline{\mathscr L'_{j-1}})
c_1(\overline{\mathscr L_{j+1}})\dots c_1(\overline{\mathscr L_d}) \\
& \leq 
\delta_{j,\infty} \deg\left( c_1(L_0)\dots c_1(L_{j-1})
c_1(L_{j+1})\dots c_1(L_d) \mid Z \right). 
\end{align*}
Adding these contributions, we get one of the desired upper bound,
and the other follows by symmetry.
\end{proof}

\begin{defi}
An adelic metric on a line bundle $L$ on~$X$ is said to 
be \emph{semipositive}\index{Semipositive>metric}\index{Metric>semipositive=(Semipositive ---)} if it is a limit of a sequence of semipositive algebraic adelic metrics on~$L$.
\end{defi}

Let $\bPic{}^+(X)$\glossary{$\overline{\mathrm{Pic}}{}^+(X)$: semi-positive isometry classes} be the set of all isomorphism classes of line bundles
endowed with a semipositive metric. We endow it with the topology
for which the set of adelic metrics on a given line bundle 
is open and closed in~$\bPic{}^+(X)$ 
and endowed with the topology induced by the distance of adelic metrics.

It is a submonoid of~$\bPic(X)$;
moreover, its image in~$\Pic(X)$ consists of (isomorphism classes
of) nef line bundles on~$X$. I thank the referee for pointing
out an example \citep[1.7]{demailly-peternell-schneider1994}
of a nef line bundle on a complex projective variety admitting 
no smooth semipositive metric, as well as for the numerous
weaknesses of exposition he/she pointed out.

\begin{coro}
Let $Z$ be an integral closed subscheme of~$X$, let $d=\dim(Z)$.
The arithmetic degree map extends uniquely to a continuous 
function $ \bPic{}{}^+(X)^{d+1} \to \R$.
This extension is multilinear and symmetric.
\end{coro}
\begin{proof}
This follows from proposition~\ref{prop.continuity-height}
and from the classical extension theorem of uniformly continuous maps.
\end{proof}

\begin{defi}
Let $X$ be a projective $\Q$-scheme and let $L$ be a line bundle on~$X$.
An adelic metric on~$L$ is said to be \emph{admissible}\index{Admissible>metric}\index{Metric>admissible=(Admissible ---)}
if there exists two line bundles endowed with semipositive adelic  metrics,
$\overline {M_1}$ and~$\overline {M_2}$, such that
$\overline L\simeq \overline{M_1}\otimes \overline{M_2}{}^{-1}$.
\end{defi}

More generally, we say that a $v$-adic metric on~$L$
is admissible if it is the $v$-adic component of an admissible adelic metric on~$L$.
The set of all admissible adelically metrized line bundles on~$X$ is denoted
by~$\bPic{}^{\adm}(X)$; it is the subgroup generated
by $\bPic{}^+(X)$.\glossary{$\overline{\mathrm{Pic}}{}^{\mathrm{adm}}(X)$: Admissible classes} 

By construction, the arithmetic intersection product
extends by linearity to $\bPic{}^\adm(X)$.
We use the notation $\hdeg(\hc_1(\overline L_0)\cdots\hc_1(\overline L_d)\mid Z)$ for the arithmetic degree 
of a $d$-dimensional integral closed subscheme~$Z$ of~$X$
with respect to admissible
adelically metrized line bundles $\overline L_0,\dots,\overline L_d$.

This gives rise to a natural notion of height parallel to
that given in definition~\ref{defi.height}.

\begin{exem}
Let us retain the context and notation of example~\ref{exem.algebraic-dynamics}.
Let us moreover assume that $L$ is ample and let us prove that
the canonical adelic metric on~$L$ is semipositive.

We make the observation that if 
$\norm\cdot $ is an algebraic adelic metric on~$L$  induced by a relatively
semipositive hermitian line bundle $\overline{\mathscr L}$ on a proper flat model~$\mathscr X$ of~$X$, then the metric $\norm\cdot^f$ is again relatively
semipositive. Indeed, the normalization of~$\mathscr X$
in the morphism $f\colon X\to X$ furnishes a proper flat scheme~$\mathscr X'$ 
over~$\Z$ such that $\mathscr X'_\Q=X$ and a morphism $\phi\colon \mathscr X'\to\mathscr X$ that extends~$f$.
Then $\phi^*\overline{\mathscr L}$ is a relatively semipositive
hermitian line bundle on~$\mathscr X'$,
model of $L^d$,  which induces the 
algebraic adelic metric $\norm\cdot^f$ on~$L$.

Starting from a given algebraic adelic metric induced by a relatively
semipositive model (for example, a relatively ample one), the proof of Picard's theorem invoked in example~\ref{exem.algebraic-dynamics}
proves that the sequence of adelic metrics obtained by the iteration
of the operator $\norm\cdot\mapsto \norm\cdot^f$ converges to the
unique fixed point.
Since this iteration preserves algebraic adelic metrics 
induced by a relatively semipositive model,
the canonical adelic metric on~$L$ is semipositive, as claimed.

For a generalization of this construction, 
see theorem~4.9 of~\cite{yuan-zhang2017}.
\end{exem}

\section{Measures}

\begin{defi}
Let $X$ be a projective $\Q$-scheme.
A function $\phi\in \mathscr C(X_\ad,\R)$ is said to be admissible\index{Admissible>function}
if the adelically metrized  line bundle $\mathscr O_X(\phi)$
is admissible.
\end{defi}

The set $\mathscr C_\adm(X_\ad,\R)$
of admissible functions $(\phi_v)$ is a real vector subspace
of $ \mathscr C_\cpct(X_\ad,\R)$. One has an exact sequence
\begin{equation}
 \Gamma(X,\mathscr O_X^\times) \to \mathscr C_\adm(X_\ad,\R) \to \bPic{}^\adm(X) \to \Pic(X) \to0
\end{equation}
analogous to~\eqref{eq.pic}.

More generally, we say that 
a function $\phi_v\in\mathscr C(X_v^\an,\R)$
is admissible if it is the $v$-adic component 
of an admissible function $\phi=(\phi_v)$. This defines
a real vector subspace~$\mathscr C_\adm(X_v^\an,\R)$
of~$\mathscr C(X_v^\an,\R)$.

\begin{prop}[{\citealp[theorem~7.12]{gubler1998}}]
For every place~$v\in S$, the subspace $\mathscr C_\adm(X_v^\an,\R)$
is dense in~$\mathscr C(X_v^\an,\R)$.

The space  $\mathscr C_\adm(X_\ad,\R)$ of admissible functions
is dense in 
$\mathscr C_\cpct(X_\ad,\R)$.
\end{prop}
\begin{proof}
Observe that $X_v^\an$ is a compact topological space. 
By corollary~7.7 and lemma~7.8 of \cite{gubler1998}, 
the subspace of $\mathscr C_{\adm}(X_v^\an,\R)$ corresponding
to algebraic $v$-adic metrics on~$L$ separates
points and is stable under~$\sup$ and~$\inf$. 
The first part of the proposition thus follows from Stone's density theorem.

The second part follows from the first one and a straightforward argument.
\end{proof}

\begin{theo}
Let $Z$ be an integral closed subscheme of~$X$,
let $d=\dim(Z)$, let $\overline{L_1},\dots,\overline{L_d}$
be admissible adelically  metrized line bundles on~$X$.

\begin{enumerate}
\item
There exists  a unique measure 
$c_1(\overline{L_1})\dots c_1(\overline{L_d})\delta_Z$
on~$X_\ad$ such that
\[ \int_{X_\ad} \phi_0 c_1(\overline{L_1})\dots c_1(\overline{L_d})\delta_Z
= \hdeg (\hc_1( \mathscr O_X(\phi_0)) \hc_1(\overline{L_1})\dots \hc_1(\overline{L_d})\mid Z
)
 \]
for every compactly supported
admissible function~$\phi_0$ on~$X_\ad$.

\item
This measure is supported on~$Z_\ad$; for every place~$v$ of~$S$,
the mass of its restriction to~$X_v^\an$ is equal to
\[ 
 \int_{X_v^\an}  c_1(\overline{L_1})\dots c_1(\overline{L_d})\delta_Z
= \deg (c_1(L_1)\cdots c_1(L_d)\mid Z). \]
If $\overline{L_1},\dots,\overline{L_d}$ are semipositive, then this
measure is positive.

\item
The induced map from~$\bPic^\adm(X)^d$ to the space
$\mathscr M(X_\ad)$ of measures on~$X_\ad$
is $d$-linear and symmetric.

\item
Every admissible function is integrable for this measure.
\end{enumerate}
\end{theo}
\begin{proof}
Let us first assume that $\overline{L_1},\dots,\overline{L_d}$
are semipositive.
It then follows from the definition of the arithmetic intersection
degrees that the map
\[ \phi_0\mapsto 
\hdeg (\hc_1( \mathscr O_X(\phi_0)) \hc_1(\overline{L_1})\dots \hc_1(\overline{L_d})\mid Z
)\]
is a positive  linear form on~$\mathscr C_\adm(X_\ad,\R)$.
By the density theorem, it extends uniquely to a positive
linear form on~$\mathscr C_c(X_\ad,\R)$, which then corresponds to
an inner regular, locally finite, positive Borel measure on~$X_\ad$.

The rest of the theorem follows from this.
\end{proof}

\begin{rema}
\begin{enumerate}
\item
At archimedean places, the construction of the measure
$c_1(\overline{L_1})\dots c_1(\overline{L_d})\delta_Z$
shows that it coincides with the measure defined
by~\cite{bedford-t82} and~\cite{demailly1985}.

\item
At finite places, it has been first given in~\cite{chambert-loir2006}.
By approximation, the definition of the measure
in the case of a general semipositive $p$-adic metric is then deduced
from the case of algebraic metrics,
given by a model~$(\mathscr X,\mathscr L)$.
In this case, the measure
$c_1(\overline{L_1})\dots c_1(\overline{L_d})\delta_Z$
on~$X_p^\an$ has finite support. Let us describe it when
 $Z=X$ and the model~$\mathscr X$ (the general case follows).
For each component~$\mathscr Y$ of~$\mathscr X\otimes\F_p$, there exists
a unique point~$y\in X_p^\an$ whose specialization is the
generic point of~$\mathscr Y$. The contribution of the point~$y$
to the measure is then equal to
\[ m_{\mathscr Y} \deg(
c_1(\mathscr L_1)\dots c_1(\mathscr L_d)\mid \mathscr Y ), \]
where $m_{\mathscr Y}$ is the multiplicity of~$\mathscr Y$
in the special fiber, that is, the length of the ideal~$(p)$ at the
generic point of~$\mathscr Y$.
\end{enumerate}
\end{rema}

\begin{exem}
Let $X$ be an abelian variety of dimension~$d$ over a number field~$F$.
Let $\overline L$ be an ample line bundle equipped with a canonical
adelic metric; let us then describe the measure $c_1(\overline L)^d$
on~$X_v^\an$, for every place~$v\in S$.
For simplicity, we assume that $F=\Q$.

\begin{enumerate}
\item
First assume $v=\infty$. Then $X_\infty^\an$ is the quotient,
under complex conjugation, of the complex torus~$X(\C)$,
and the canonical measure on~$X_\infty^\an$ is the direct image
of the unique Haar measure on~$X(\C)$ with total mass
$\deg(c_1(L)^d\mid X) $.

\item
The situation is more interesting in the case of a finite place~$p$.

If $X$ has good reduction at~$p$, that is, if it extends to an abelian
scheme~$\mathscr X$ over~$\Z_p$, then the canonical measure is supported at
the unique point of~$X_p^\an$ whose specialization is the generic
point of~$\mathscr X\otimes\F_p$.

Let us assume, on the contrary, 
that $X$ has (split) totally degenerate reduction.
In this case, the uniformization theory of abelian varieties
shows that $X_p^\an$ is the quotient of a torus~$(\gm^d)^\an$
by a lattice~$\Lambda$. The definition of $(\gm^d)^\an$ shows that this
analytic space contains a canonical $d$-dimensional real vector space~$V$,
and $V/\Lambda$ is a real $d$-dimensional
torus~$S(X_p^\an)$ contained in~$X_p^\an$,
sometimes called its skeleton. 
\cite{gubler2007b} has shown that the measure~$c_1(\overline L)^d$
on~$X_p^\an$ coincides with the Haar measure on~$S(X_p^\an)$
with total mass $\deg(c_1(L)^d\mid X)$.

The general case is a combination of these two cases,
see~\cite{gubler2010}.
\end{enumerate}
\end{exem}

\begin{rema}
At finite places, the theory described in this section
defines measures $c_1(\overline L_1)\dots c_1(\overline L_d)\delta_Z$
without defining the individual components~$c_1(\overline L_1)$,\dots, $c_1(\overline L_d)$, and $\delta_Z$.  
In \cite{chambert-loir-ducros2012}, we propose a theory
of real differential forms and currents on Berkovich analytic spaces  
that allows a more satisfactory analogy with the theory on complex spaces.
In particular, we provide an analogue of the Poincaré--Lelong equation,
and semipositive metrized line bundles possess a curvature current
$c_1(\overline L)$
(curvature form in the ``smooth'' case) whose product can be defined
and coincides with the measures of
the form $c_1(\overline L_1)\dots c_1(\overline L_d)\delta_Z$
that we have discussed.
\end{rema}

\section{Volumes}

\subsection{}
Let $X$ be a proper $\Q$-scheme and 
let $\overline L$ be a line bundle endowed with an adelic metric.

The Riemann-Roch space $\mathrm H^0(X,L)$ 
is a finite dimensional $\Q$-vector space.
For every place $v\in S$, we endow it with a $v$-adic semi-norm:
\[ \norm s_v = \sup_{x\in X_v^\an} \norm{s(x)} \]
for $s\in \mathrm H^0(X,L)$. If $X$ is reduced, then this is a norm;
let then $B_v$ be its unit ball.

In the case where the metric on~$L$ is algebraic, given by a model~$\mathscr L$
on a \emph{normal} model~$\mathscr X$ of~$X$, 
the real vector space $\mathrm H^0(X,L)_\R$ admits
a natural lattice $\mathrm H^0(\mathscr X,\mathscr L)$
and the norm~$\norm\cdot_\infty$ at the archimedean place.
The real space $\mathrm H^0(X,L)_\R$ can then be endowed
with the unique Haar measure for which~$B_\infty$ has volume~$1$,
and we can then consider the covolume of the lattice 
$\mathrm H^0(\mathscr X,\mathscr L)$. Minkowski's first
theorem will then eventually provide non-trivial
elements of $\mathrm H^0(\mathscr X,\mathscr L)$ of
controlled norm at the archimedean place.

The following construction, due to~\cite{bombieri-v1983},
provides the analogue in the adelic setting.

Let $\mathbf A$ be the ring of adeles of~$\Q$ and let $\mu$
be a Haar measure on $\mathrm H^0(X,L)\otimes\A$. 
Then $\prod_{v\in S}B_v$  has finite strictly positive volume in $\mathrm H^0(X,L)_\A$,
and one defines
\begin{equation}
\chi( X,\overline L) =  - \log 
\left( \frac{\mu(\mathrm H^0(X,L)\otimes\A/\mathrm  H^0(X,L))}{\mu(\prod_v B_v)}
\right).
\end{equation}
This does not depend on the  choice of the Haar measure~$\mu$.

One also defines
\begin{equation} \widehat {\mathrm H}^0(X,\overline L)= \{s\in \mathrm H^0(X,L) \sozat
   \norm s_v\leq 1 \text{ for all $v\in S$}\}. \end{equation}
This is a finite set.
We then let
\begin{equation} \widehat {\mathrm h}^0(X,\overline L) = \log \left( \Card (\widehat{\mathrm H}^0(X,\overline L))\right). \end{equation}

The following inequality is an adelic analogue of Blichfeldt's theorem,
itself a close companion to the adelic analogue of Minkowski's first theorem 
proved by~\cite{bombieri-v1983}.
We refer to corollaire~2.12 of~\cite{Gaudron-2009} for details. 
\begin{lemm}
One has 
\[ \chi (X,\overline L) \leq \widehat {\mathrm h}^0(X,\overline L)
+ \mathrm h^0(X,L)\log(2). \]
\end{lemm}

\subsection{}
The volume and the $\chi$-volume of~$\overline L$ are defined by the formulas\index{Volume>of a metrized line bundle}:
\glossary{$\widehat{\text{Vol}}(X,\overline L)$: Volume of a metrized line bundle}
\glossary{$\widehat{\text{Vol}}_\chi(X,\overline L)$: $\chi$-Volume of a metrized line bundle}
\begin{gather}
\hvol(X,\overline L) = \limsup_{n\to\infty} \frac{\widehat{\mathrm h}^0(X,\overline L^n)}{n^{d+1}/(d+1)!} \\
\hvolchi(X,\overline L) = \limsup_{n\to\infty} \frac{\chi(X,\overline L^n)}{n^{d+1}/(d+1)!} .
\end{gather}
One thus has the inequality
\begin{equation}
 \hvolchi(X,\overline L)\leq \hvol(X,\overline L).
\end{equation}

In fact, it has been independently shown by~\cite{yuan2009}
and~\cite{chen2010b}
that the volume is in fact a limit.

The relation between volumes and heights follows from the following
result.
\begin{lemm}\label{lemm.volume-height}
Assume that $L$ is ample. Then,
for every real number~$t$
such that 
\[ t < \dfrac{\hvolchi(X,\overline L)}{(d+1)\vol(X,L)}, \]
the set of closed points $x\in X$ such that $h_{\overline L}(x)\leq t$
is not dense for the Zariski topology.
\end{lemm}
\begin{proof}
Consider the adelically metrized line bundle~$\overline L(-t)$, 
whose metric at the archimedean place has been multiplied by~$e^{t}$.
It follows from the definition of the $\chi$-volume that 
\[ \hvol\nolimits_\chi (X,\overline L(-t))=\hvol\nolimits_\chi (X,\overline L) - (d+1) t \vol(X,L). \]
Indeed, for every finite place~$p$,
changing~$\overline L$ to~$\overline L(-t)$ does not modify the balls~$B_p$ in $\mathrm H^0(X,\overline L^n)\otimes_\Q\Q_p$,
while it dilates it by the ratio~$e^{-nt}$ at the archimedean place, so that its volume
is multiplied by $e^{-nt \dim(\mathrm H^0(X,L^n))}$.

Consequently,
\[   \hvol (X,\overline L(-t))\geq \hvolchi(X,\overline L(-t))
\geq \hvolchi (X,\overline L) - (d+1) t \vol(X,L) >0 . \]
In particular, there exists an integer~$n\geq 1$ and a nonzero section $s\in \mathrm H^0(X,\overline L^n)$
such that $\norm s_p\leq 1$ for all finite places~$p$, and $\norm s_\infty\leq e^{-nt}$.
Let now $x\in X$ be a closed point that is not contained in $\abs{\div(s)}$; 
one then has
\[ h_{\overline L}(x) = \sum_{v\in S} \int_{X_v^\an} \log \norm{s}_v^{-1/n}\delta_v(x) \geq  t, \]
whence the lemma.
\end{proof}

\begin{theo}\label{theo.hilbert-samuel}
Assume that $\overline L$ is semipositive.
Then one has 
\begin{equation}
\hvol(X,\overline L)=\hvol\nolimits_\chi(X,\overline L) = \hdeg \left(
   \hc_1(\overline L)^{d+1}\mid X\right).
\end{equation}
\end{theo}
This is the arithmetic Hilbert--Samuel formula,
due to~\cite{gillet-s88,bismut-vasserot1989}
when $X$ is smooth 
and the adelic metric of~$\overline L$ is algebraic.
\cite{abbes-bouche:1995} later gave an alternative proof.
In the given generality, the formula 
is a theorem of~\cite{zhang95,zhang95b}.

\begin{theo}
\begin{enumerate}
\item
The function $\overline L \mapsto \hvol(X,\overline L)$
extends uniquely 
to a continuous function on the real vector space~$\bPic(X)\otimes_\Q\R$.

\item
If $\hvol(X,\overline L)>0$, then  $\hvol$ is
differentiable at~$\overline L$.

\item
If $\overline L$ is semipositive, then $\hvol$ and $\hvol\nolimits_\chi$
are differentiable at $\overline L$, with differential
\[ \overline M \mapsto (d+1) \hdeg (\hc_1 (\overline L)^d \hc_1(\overline M)\mid X).\]
\end{enumerate}
\end{theo}

This theorem is proved by~\cite{chen2011} as a consequence
of results of~\cite{yuan2008,yuan2009}.
It essentially reduces
from the  preceding one in the case $\overline L$ is defined by an ample
line bundle on a model of~$X$, and its metric has strictly
positive curvature. 
Reaching the ``boundary'' of the cone of semipositive admissible
metrized line bundles was the main result of~\cite{yuan2008}
who proved that for every admissible metrized line bundle $\overline M$
and every large enough integer~$t$, one has
\[ 
 t^{-(d+1)} \hvolchi(X,\overline L^t\otimes \overline M)
   \geq \hvolchi(X,\overline L)
 + \frac1t (d+1) \hdeg (\hc_1(\overline L)^d \hc_1(\overline M)\mid X) 
 + \mathrm o(1/t).
\]
It is this inequality, an arithmetic analogue of an inequality of Siu,
will be crucial for the applications
to equidistribution in the next section.

\section{Equidistribution}

The main result of this section is the equidistribution 
theorem~\ref{theo.equidistribution}.
It has been first proved in the case $v=\C$ by~\cite{szpiro-u-z97},
under the assumption that the given archimedean metric is smooth and has a strictly positive curvature form,  and the general case is due to~\cite{yuan2008}.
However, our presentation derives it 
from a seemingly more general result, lemma~\ref{lemm.equidistribution},
whose proof, anyway,  closely follows their methods. 
Note that for the application to Bogomolov's conjecture 
in~\S\ref{sec.bogomolov},
the initial theorem of~\cite{szpiro-u-z97} is sufficient.

\begin{defi}
Let $X$ be a proper $\Q$-scheme, let $\overline L$ be an ample line
bundle on~$X$ endowed with an admissible
adelically metric.
Let $(x_n)$ be a sequence of closed points of~$X$.

\begin{enumerate}
\item One says that $(x_n)$ is \emph{generic}\index{Generic sequence}\index{Sequence>generic=(Generic ---)} if 
for every strict closed 
subscheme~$Z$ of~$X$, the set of all $n\in\N$ such that $x_n\in Z$ is finite;
in other words, this sequence converges to
the generic point of~$X$. 
\item One says that $(x_n)$ is \emph{small}\index{Small sequence}\index{Sequence>small=(Small ---)} (relative to~$\overline L$) if
\[ h_{\overline L}(x_n) \to h_{\overline L}(X). \]
\end{enumerate}
\end{defi}

\begin{lemm}\label{lemm.equidistribution}
Let $X$ be a proper $\Q$-scheme, let $d=\dim(X)$, let $\overline L$
be a semipositive adelically metrized line bundle on~$X$
such that $L$ is ample.
Let $(x_n)$ be a generic sequence of closed points of~$X$
which is small relative to~$\overline L$.
For every line bundle~$\overline M$ on~$X$ endowed with an admissible adelic metric, one has
\[
\lim_{n\to\infty} h_{\overline M} (x_n)  
= \frac{\hdeg(\hc_1(\overline L)^d\hc_1(\overline M)\mid X)}{\deg_L(X)}
- {d} h_{\overline L}(X) \frac{\deg(c_1(L)^{d-1}c_1(M)\mid X)}{\deg_L(X)}. 
\]
\end{lemm}
\begin{proof}
Since $L$ is ample, $L^t\otimes M$ is ample for every large enough integer~$t$,
and the classical Hilbert-Samuel formula implies that
\[ 
\frac1{t^d} \vol(X,L^t\otimes M)  = \deg (c_1(L)^d\mid X) 
+ d t^{-1} \deg( c_1(L)^{d-1}c_1(M)\mid X) + \mathrm O(t^{-2})
\]
when $t\to\infty$.
Since $\overline L$ is semipositive and $L$ is ample, 
the main inequality of~\cite{yuan2008} implies that
\begin{multline*}
 \frac1{t^{d+1}} \hvolchi (X,\overline L^t \otimes \overline M)
  \geq \hdeg(\hc_1(\overline L)^{d+1}\mid X) 
\\
{} + (d+1) t^{-1} \hdeg (\hc_1(\overline L)^d\hc_1(\overline M)\mid X) 
 + \mathrm o(t^{-1}) .
\end{multline*}
Consequently, when $t\to \infty$, one has 
\begin{multline*}
\frac{\hvolchi(X,\overline L^t\otimes M)}{\vol(X,L^t\otimes M)} 
\geq t \frac{\hdeg(\hc_1(\overline L)^{d+1}\mid X)}{\deg(c_1(L)^d\mid X)} \\
{} +  (d+1) \frac{ \hdeg(\hc_1(\overline L)^d\hc_1(\overline M)\mid X)}{\deg(c_1(L)^d\mid X)}  \\
{} 
 -d \frac{ \hdeg(\hc_1(\overline L)^{d+1}\mid X)}{\deg (c_1(L)^d\mid X)}
\frac { \deg(c_1(L)^{d-1}c_1(M)\mid X)}{\deg(c_1(L)^d\mid X)}
+ \mathrm o(1).
\end{multline*}

The sequence~$(x_n)$ is generic, hence lemma~\ref{lemm.volume-height}
furnishes the inequality:
\[ \liminf_n h_{\overline L^t\otimes\overline M}(x_n)
\geq  \frac{\hvolchi(X,\overline L^t\otimes M)}{(d+1)\vol(X,L^t\otimes M)}.
\]
We observe that
\[  \liminf_n h_{\overline L^t\otimes\overline M} (x_n)
= t \lim  h_{\overline L}(x_n ) + \liminf_n h_{\overline M}(x_n), \]
so that, when $t\to\infty$, we have
\begin{multline*}
 \liminf_n h_{\overline M}(x_n) \geq 
 \frac{ \hdeg(\hc_1(\overline L)^d\hc_1(\overline M)\mid X)}{\deg(c_1(L)^d\mid X)} \\
{}
- \frac d{d+1} \frac{ \hdeg(\hc_1(\overline L)^{d+1}\mid X)}{\deg (c_1(L)^d\mid X)}
\frac { \deg(c_1(L)^{d-1}c_1(M)\mid X)}{\deg(c_1(L)^d\mid X)}.
\end{multline*}

Applying this inequality for $\overline M^{-1}$
shows that $\limsup_n h_{\overline M}(x_n)$ is bounded above
by its right hand side. The lemma follows.
\end{proof}

\subsection{}
Let $X$ be a proper $\Q$-scheme. Let $v\in S$ be a place of~$\Q$.

Let $x\in X$ be a closed point.
Let $F=\kappa(x)$; this is a finite extension of~$\Q$,
and there are exactly $[F:\Q]$ geometric points on $X(\C_v)$
whose image is~$x$, permuted by the Galois group $\Gal(\C_v/\Q_v)$.
We consider the corresponding ``probability measure'' in $X(\C_v)$,
giving mass $1/[F:\Q]$ to each of these geometric points,
and let $\delta_v(x)$ be its image under the natural map $X(\C_p)\to X_v^\an$.

By construction, $\delta_v(x)$ 
is a probability measure on~$X_v^\an$ with finite support,  
a point of~$X_v^\an$
being counted proportionaly to the number of its liftings to a geometric point
supported by~$x$.

The space of measures on~$X_v^\an$ is the dual of the space
$\mathscr C(X_v^\an,\R)$; it is  endowed with the topology
of pointwise convergence (``weak-* topology'').

\begin{theo}\label{theo.equidistribution}
Let $X$ be a proper $\Q$-scheme, let $d=\dim(X)$, let~$\overline L$
be a semipositive adelically metrized line bundle on~$X$
such that $L$ is ample.
Let $(x_n)$ be a generic sequence of closed points of~$X$
which is small relative to~$\overline L$.
Then for each place $v\in S$, the sequence of measures $(\delta_v(x_n))$
on~$X_v^\an$ converges to the measure  $c_1(\overline L)^d/\deg(c_1(L)^d\mid X)$.
\end{theo}
%
\begin{proof}
Let $\mu_{\overline L}$ denote the probability measure
$ c_1(\overline L)^d/\deg_L(X)$ on~$X_v^\an$ and
let $f\in\mathscr C(X_v^\an,\R)$ be an admissible function,
extended by zero to an element of $\mathscr C_\adm(X_\ad,\R)$.
We apply lemma~\ref{lemm.equidistribution} to the metrized line bundle~$\overline M=\mathscr O_X(f)$ whose underlying line bundle on~$X$ is trivial.
For every closed point~$x\in X$, one has 
\[ h_{\overline M}(x) = \int_{X_v^\an} f \delta_v(x). \]
Moreover,
\[ \hdeg (\hc_1(\overline L)^d \hc_1(\overline M)\mid X)
 =  \int_{X_v^\an} f  c_1(\overline L)^d . \]
It thus follows from lemma~\ref{lemm.equidistribution} that
\begin{align*}
 \lim_{n\to\infty} \int_{X_v^\an} f \delta_v(x_n)
&  = \frac{\hdeg(\hc_1(\overline L)^d\hc_1(\overline M)\mid X)}{\deg(c_1(L)^d\mid X)}\\
& = \frac1{\deg(c_1(L)^d\mid X)} \int_{X_v^\an} f  c_1(\overline L)^d . 
\end{align*}
The case of an arbitrary continuous  function on~$X_v^\an$
follows by density.
\end{proof}

\section{The Bogomolov conjecture}\label{sec.bogomolov}

\subsection{}
Let $X$ be an abelian variety over a number field~$F$,
that is, a projective connected algebraic variety over~$F$
which is endowed with an additional structure of an algebriac group.
For every integer~$n$, we write $[n]$ for the multiplication by~$n$-morphism 
on~$X$.

Let us first explain how the theory of canonical adelic metrics
allows to extend the Néron--Tate height to arbitrary
integral closed subschemes. For alternative and independent presentations,
see~\cite{philippon91}, \cite{gubler1994}, \cite{bost-g-s94}.

Let $L$ be a line bundle on~$X$ trivialized at the origin.

If $L$ is even ($[-1]^*L\simeq L$), the
theorem of the cube implies that there exists,
for every integer~$n$, a unique isomorphism
$[n]^*L\simeq L^{n^2}$ which is compatible with the trivialization
at the origin.
If $n\geq 2$, then by example~\ref{exem.algebraic-dynamics},
it admits a unique adelic metric
for which this isomorphism is an isometry.
is an isometry; these isomorphisms are then all isometries.

Similarly, 
if $L$ is odd ($[-1]^*L\simeq L^{-1})$, 
then it admits a unique adelic metric 
for which the canonical isomorphisms
$[n]^*L\simeq L^{n}$
are isometries, for all integer~$n$.

In general, one can write $L^2\simeq (L\otimes [-1]^*L)\otimes (L\otimes [-1]^*L^{-1})$, as the sum of an even and an odd line bundle,
and this endows~$L$ with an adelic metric.
This adelic metric is called the \emph{canonical adelic metric}\index{Canonical>adelic metric} on~$L$
(compatible with the given trivialization at the origin).

If $L$ is ample and even, then the canonical adelic metric on~$L$
is semipositive. This implies that the canonical adelic metric
of an arbitrary even line bundle is admissible.

Assume that $L$ is odd. Fix an even ample line bundle~$M$.
Up to extending the scalars,
there exists a point~$a\in X(F)$ such that $L\simeq \tau_a^*M\otimes M^{-1}$,
where $\tau_a$ is the translation by~$a$ on~$X$.
Then there exists a unique isomorphism $L\simeq \tau_a^*M \otimes M^{-1}
\otimes M_a^{-1}$ which is compatible with the trivialization at the origin,
and this gives rise to an isometry $\overline L\simeq \tau_a^*\overline M\otimes \overline M^{-1}\otimes\overline M_a^{-1}$. In particular,
the adelic metric of~$\overline M$ is admissible.
In fact, it follows from a construction of~{Künnemann}  that
it is even semipositive, see~\cite{chambert-loir99}.

\subsection{}
In particular, let us consider an ample even line bundle~$L$ on~$X$
endowed with a canonical adelic metric.
This furnishes a \emph{height}
\[ h_{\overline L}(Z) = \frac{\hdeg(\hc_1(\overline L)^{d+1}\mid Z)}{(d+1)\deg(c_1(L)^d\mid Z)}, \]
for every integral closed subscheme~$Z$ of~$X$, where $d=\dim(Z)$.

In fact, if $(\mathscr X,\overline{\mathscr L})$ is any 
model of $(X,L)$, one has
\[ h_{\overline L}(Z) = \lim_{n\to+\infty} n^{-2} h_{\overline{\mathscr L}}([n](Z)), \]
which shows the relation of  the point of view of adelic metrics 
with Tate's definition of the Néron--Tate height,
initially defined on closed points.
This formula also implies that $h_{\overline L}$ is positive.

More generally, if $Z$ is an integral closed subscheme of~$X_{\overline F}$,
we let $h_{\overline L}(Z)=h_{\overline L}([Z])$, where $[Z]$
is its Zariski-closure in~$X$ (more precisely, 
the smallest closed subscheme of~$X$ such that
$[Z]_{\overline F}$ contains~$Z$).

\begin{lemm}
The induced height function $h_{\overline L}\colon X(\overline F)\to \R$ 
is a positive quadratic form. It induces a positive definite quadratic
form on $X(\overline F)\otimes\R$.
In particular, a point $\in X(\overline F)$ satisfies $h_{\overline L}(x)=0$
if and only if $x$ is a torsion point.
\end{lemm}
\begin{proof}
For $I\subset\{1,2,3\}$, let $p_I\colon X^3\to X$ be the morphism
given by $p_I(x_1,x_2,x_3)=\sum_{i\in I}x_i$.
The cube theorem asserts that the line bundle
\[ \mathscr D_3(L)=\bigotimes_{\emptyset\neq I\subset\{1,2,3\}} (p_I^* L)^{(-1)^{\Card(I)-1}}
\]
on~$X^3$ is  trivial, and admits a canonical trivialisation.
The adelic metric of~$\overline L$ endows
it with an adelic metric which satisfies $[2]^*\mathscr D_3(\overline L)
\simeq \mathscr D_3(\overline L)^4$, hence is the trivial metric.
This implies the following relation on heights:
\[ h_{\overline L}(x+y+z) - h_{\overline L}(y+z)-h_{\overline L}(x+z)-h_{\overline L}(x+y)+h_{\overline L}(x)+h_{\overline L}(y)+h_{\overline L}(z)\equiv 0\]
on $X(\overline F)^3$.
Consequently, 
\[ (x,y) \mapsto h_{\overline L}(x+y)-h_{\overline L}(x)-h_{\overline L}(y)\]
is a symmetric bilinear form on~$X(\overline F)$.
Using furthermore that $h_{\overline L}(-x)=h_{\overline L}(x)$
for all~$x$, 
we deduce that $h_{\overline L}$  is a quadratic form on~$X(\overline F)$.

Since $L$ is ample, $h_{\overline L}$ is bounded from below.
The formula $ h_{\overline L} (x)= h_{\overline L}(2x)/4$ then implies
that $h_{\overline L}$ is positive.
By what precedes, it induces a  positive  quadratic form on~$X(\overline F)_\R$.

Let us prove that it is in fact positive definite.
By definition, it suffices
that its restriction to  the subspace generated
by finitely many points $x_1,\dots,x_m\in X(\overline F)$ is positive definite.
Let $E$ be a finite extension of~$F$ such that $x_1,\dots,x_m\in X(E)$.
On the other hand, Northcott's theorem
implies that for every real number~$t$, the set of $(a_1,\dots,a_m)\in\Z^m$
such that $h_{\overline L}(a_1x_1+\dots+a_mx_m)\leq t$ is finite. 
One deduces from that the asserted positive definiteness.
\end{proof}

\begin{defi}
A \emph{torsion subvariety}\index{Torsion subvariety}\index{Subvariety>torsion=(Torsion ---)} of~$X_{\overline F}$
is a subvariety of the form $a+Y$, where $a\in X(\overline F)$ 
is a torsion point and $Y$ is an abelian subvariety of~$X_{\overline F}$.
\end{defi}

\begin{theo}\phantomsection\label{theo.bogomolov}
\begin{enumerate}\def\theenumi{\alph{enumi}}\def\labelenumi{\theenumi)}
\item 
Let $Z$ be an integral closed subscheme of~$X_{\overline F}$.
One has $h_{\overline L}(Z)=0$ if and only if $Z$ is a torsion
subvariety of~$X_{\overline F}$.

\item
Let $Z$ be an integral closed subscheme of~$X_{\overline F}$
which is not a torsion subvariety.
There exists a strictly positive real number~$\delta$
such that the set
\[ \{ x\in Z(\overline F)\sozat h_{\overline L} (x) \leq\delta \} \]
is not Zariski-dense in~$Z_{\overline F}$.
\end{enumerate}
\end{theo}

Assertion~\emph a) has been independently conjectured by
\cite{philippon91,philippon95} and~\cite{zhang95b}.
Assertion~\emph b) has been conjectured by~\cite{bogomolov80b}
in the particular case where $Z$ is a curve of genus~$g\geq 2$
embedded in its jacobian variety; for this reason, it is called
the ``generalized Bogomolov conjecture''.
The equivalence of~\emph a) and~\emph b) is a theorem of~\cite{zhang95b}.
In fact, the implication \emph a)$\Rightarrow$\emph b)
already follows from theorem~\ref{theo.hilbert-samuel}
and lemma~\ref{lemm.volume-height}. Assume indeed that
\emph a) holds and let $Z$ be an integral closed
subscheme of~$X$ which is not a torsion subvariety; 
by \emph a), one has $h_{\overline L}(Z)>0$,
and one may take for~$\delta$ any real number such that
$0<\delta< h_{\overline L}(Z)$.

Theorem~\ref{theo.bogomolov} has been proved by~\cite{zhang98},
following a breakthrough of~\cite{ullmo98} who treated the case
of a curve embedded in its jacobian; their proof makes use
of the equidistribution theorem.
Soon after, \cite{david-p98} gave an alternative proof;
when $Z$ is not a translate of an abelian subvariety, their
proof provides a strictly positive lower bound 
for $h_{\overline L}(Z)$ (in~\emph a))
as well as an explicit real number~$\delta$ (in~\emph b))
which only depends on the dimension and the degree of~$Z$ with respect to~$L$.

As a corollary of theorem~\ref{theo.bogomolov}, one obtains a new
proof of the Manin--Mumford conjecture in characteristic zero,
initially proved by~\cite{raynaud83c}.
\begin{coro}
Let $X$ be an abelian variety over an algebraically closed field  of characteristic zero, let $Z$ be an integral closed subscheme of~$X$
which is not a torsion subvariety. Then the set of torsion points of~$X$
which are contained in~$Z$ is not Zariski-dense in~$Z$.
\end{coro}
\begin{proof}
A specialization argument reduces to the case where~$X$
is defined over a number field~$F$. In this case,
the torsion points of~$X$ are defined over~$\overline F$
and are characterized by the vanishing of their Néron--Tate height
relative to an(y) ample line bundle~$L$ on~$X$.
It is thus clear that the corollary follows from
theorem~\ref{theo.bogomolov}, \emph b).
\end{proof}

\subsection{}
The proof of theorem~\ref{theo.bogomolov}, \emph b), 
begins with the observation that the statement does not depend 
on the choice of the ample line~$L$ on~$X$. More precisely,
if $\overline M$ is another symmetric ample line bundle on~$X$ endowed with
a canonical metric, then there exists an integer~$a\geq 1$ such
that $\overline L^a\otimes\overline M^{-1}$ is ample, as well
as $\overline M^a\otimes\overline L^{-1}$. Consequently,
$h_{\overline L} \geq a^{-1} h_{\overline M}$ 
and $h_{\overline M}\geq a^{-1} h_{\overline L}$. From these
two inequalities, one deduces readily that the statement holds for~$\overline L$
if and only if it holds for~$\overline M$.

For a similar reason, if $f\colon X'\to X$ is an isogeny of abelian varieties,
then the statements for~$X$ and~$X'$ are equivalent. Let
indeed~$Z$ be an
integral closed subvariety of~$X_{\overline F}$ 
and let $Z'$ be an irreducible component of~$f^{-1}(Z)$.
Then $Z$ is a torsion subvariety of~$X_{\overline F}$
if and only if $Z'$ is a torsion subvariety of~$X'_{\overline F}$.
On the other hand, the relation $h_{f^*\overline L}(x)=h_{\overline L}(f(x))$
shows that $h_{f^*\overline L}$ has a strictly positive lower bound on~$Z'$
outside of a strict closed subset~$E'$ if and only if
$h_{\overline L}$ has a strictly positive lower bound on~$Z$
outside of the strict closed subset~$f(E')$.

\subsection{}
Building on that observation, one reduces the proof of the theorem
to the case where the stabilizer of~$Z$ is trivial.

Let indeed $X''$ be the neutral component of this stabilizer
and let $X'=X/X''$; this is an abelian variety.
By Poincaré's complete reducibility theorem, there exists an isogeny
$f\colon X'\times X''\to X$. This reduces us to the case where
$X=X'\times X''$ and $Z=Z'\times X''$, for some integral closed
subscheme~$Z'$ of~$X'_{\overline F}$.
We may also assume that $\overline L=\overline L'\boxtimes\overline L''$.
It it then clear that the statement for~$(X',Z')$ implies
the desired statement for~$(X,Z)$.

\begin{lemm}
Assume that $\dim(Z)>0$ and that its stabilizer is trivial.
Then, for every large enough integer~$m\geq 1$, the morphism
\[ f\colon Z^m\to X_{\overline F}^{m-1}, \qquad (x_1,\dots,x_m)\mapsto
 (x_2-x_1,\dots,x_m-x_{m-1}) \]
is birational onto its image but not finite.
\end{lemm}
\begin{proof}
Let $m$ be an integer and let $x=(x_1,\dots,x_m)$ 
be an $\overline F$-point of~$Z^m$.
Then a point $y=(y_1,\dots,y_m)\in Z(\overline F)^m$ belongs to
the same fiber as~$x$ if and only if $y_2-y_1=x_2-x_1,\dots$,
that is, if and only if, $y_1-x_1=y_2-x_2=\dots=y_m-x_m$. 
This identifies $f^{-1}(f(x))$ with the intersection
$(Z-{x_1})\cap \dots \dots\cap (Z-{x_m})$ of translates of~$Z$.
If $m$ is large enough and $x_1,\dots,x_m$ are well chosen in~$Z$,
then this intersection is equal to stabilizer of~$Z$ in~$X_{\overline F}$,
hence is reduced to a point.
In that case, the morphism~$f$ has a fiber reduced to a point,
hence it is generically injective.

On the other hand, the preimage of the origin~$(o,\dots,o)$
contains the diagonal of~$Z^m$, which has strictly positive dimension
by hypothesis.
\end{proof}

\subsection{}
For the proof of theorem~\ref{theo.bogomolov}, \emph b), we
now argue by contradiction and assume the existence of 
a generic sequence $(x_n)$ in~$Z(\overline F)$ 
such that $h_{\overline L}(x_n)\to 0$.

Having reduced, as explained above, to the case where the stabilizer
of~$Z$ is trivial, we consider an integer~$m\geq 1$
such that the morphism $f\colon Z^m\to X^{m-1}_{\overline F}$
is birational onto its image, but not finite.

Since the set of strict closed subschemes of~$Z$ is countable,
one can construct  a generic sequence $(y_n)$ in~$Z^m$ 
where $y_n$ is of the form $(x_{i_1},\dots,x_{i_m})$.
One has $h_{\overline L}(y_n)\to 0$,
where, by abuse of language, we write $h_{\overline L}$
for the height on~$X^m$ induced by the adelically
metrized line bundle $\overline L\boxtimes\dots\boxtimes\overline L$ on~$X^m$.
This implies that $h_{\overline L}(Z)=0$, hence the
sequence~$(y_n)$ is \emph{small}.

For every integer~$n$, let $z_n=f(y_n)$. By continuity of
a morphism of schemes, the sequence $(z_n)$
is generic in~$f(Z^m)$. Moreover, we deduce from the
quadratic character of the Néron--Tate height~$h_{\overline L}$
that $h_{\overline L}(z_n)\to 0$.
In particular, $h_{\overline L}(f(Z^m))=0$, and the sequence~$(z_n)$
is small.

Fix an archimedean place~$\sigma$ of~$F$.
Applied to the sequences~$(y_n)$ and~$(z_n)$,
the equidistribution theorem~\ref{theo.equidistribution}
implies the following convergences:
\begin{gather*}
\lim_{n\to\infty} \delta_\sigma(y_n) \propto
 c_1(\overline L\boxtimes\dots\boxtimes \overline L)^{m\dim(Z)}\delta_{Z^m}\\
\lim_{n\to\infty} \delta_\sigma(z_n) \propto
 c_1(\overline L\boxtimes\dots\boxtimes \overline L)^{\dim(f(Z^m))}\delta_{f(Z^m)},\end{gather*}
where, by $\propto$, I mean that both sides are proportional.
(The proportionality ratio is the degree of~$Z^m$, \resp of~$f(Z^m)$, with respect to the indicated measure.)
Since $f(y_n)=z_n$, we conclude that the measures
\[
f_* (c_1(\overline L\boxtimes\dots\boxtimes \overline L)^{m\dim(Z)}\delta_{Z^m})
\quad \text{and}\quad
 c_1(\overline L\boxtimes\dots\boxtimes \overline L)^{\dim(f(Z^m))}\delta_{f(Z^m)}\]
on~$f(Z^m)$ are proportional.

Recall that the archimedean metric of~$\overline L$ has the property
that it is smooth and that its curvature form~$c_1(\overline L)$
is a smooth strictly positive $(1,1)$-form on $X_\sigma(\C)$.
Consequently, on a dense smooth 
open subscheme of~$f(Z^m)$ above which $f$ is an isomorphism,
both measures are given by 
differential forms, which thus coincide there. We can pull back
them to~$Z^m$ by~$f$ and obtain a proportionality of differential forms
\[
 c_1(\overline L\boxtimes\dots\boxtimes \overline L)^{m\dim(Z)}
 \propto
 f^* c_1(\overline L\boxtimes\dots\boxtimes \overline L)^{m\dim (Z)}
\]
on~$Z_\sigma(\C)^m$. 
At this point, the contradiction appears: the differential form on the
left is strictly positive at every point, while the one on the right
vanishes at every point of~$Z_\sigma^m(\C)$ at which~$f$ is not smooth.
 
This concludes the proof of theorem~\ref{theo.bogomolov}.

\begin{rema}
The statement of~\ref{theo.bogomolov} can be asked in more general contexts
that allow for canonical heights. The case of toric varieties
has been proved by~\cite{zhang95}, while in that case
the equidistribution result is first due to~\cite{bilu97}.
The case of semiabelian varieties is due to~\cite{david-p2002},
by generalization of their proof for abelian varieties;
I had proved in~\cite{chambert-loir99} the equidistribution result 
for almost-split semi-abelian varieties, and the general case 
has just been announced by~\cite{kuehne2018}.

The general setting of algebraic dynamics $(X,f)$ is unclear. 
For a polarized dynamical system as in~\ref{exem.algebraic-dynamics},
the obvious and natural generalization proposed in~\cite{zhang95b} 
asserts that subvarieties of height zero are exactly those whose 
forward orbit is finite. However, 
{Ghioca} and~{Tucker} have  shown that 
it does not hold; see~\cite{ghioca-tucker-zhang2011} for a possible
rectification. 
The case of dominant endomorphisms of~$(\P^1)^n$ 
is a recent theorem of~\cite{ghioca-nguyen-ye2017}.
\end{rema}

\backmatter

\nocite{abbes1997,moriwaki2014a,soule-a-b-k92}
\bibliographystyle{mynat}
\bibliography{aclab,acl,grenoble}


\end{document}